\documentclass[11pt]{article}
\usepackage{thm-gcsppr}
\usepackage{thesis}

\usepackage{palatino}
\usepackage{setspace}
\usepackage{indentfirst}
\usepackage{hyperref}
\usepackage{graphicx}
\usepackage[all]{xy}

\newcommand{\anchor}{\rho} % anchor map
\newcommand{\DD}{\mathcal{D}}
\newcommand{\thalf}{\tfrac{1}{2}}
\DeclareMathOperator{\isom}{iso}
\newcommand{\diese}{^{\sharp}}

\begin{document}
\title{Generalized Complex Submanifolds}
\author{James Barton \\ Department of Mathematics \\ Pennsylvania State University \\ 
University Park, PA 16802 \\ \texttt{barton@math.psu.edu} 
\and Mathieu Sti\'enon \thanks{ESI Junior Research Fellow} \thanks{Research partially supported by NSF grant DMS03-06665.} \\ 
Departement Mathematik \\ E.T.H. Zurich \\  
8092 Zurich (Switzerland) \\ \texttt{stienon@math.ethz.ch} }
\date{}
\maketitle
\begin{abstract}
We introduce the notion of twisted generalized complex submanifolds and describe an
equivalent characterization in terms of Poisson-Dirac submanifolds.  Our characterization
recovers a result of Vaisman \cite{VAIS}.  An equivalent characterization is also given
in terms of spinors.  As a consequence, we show that the fixed locus of an involution
preserving a twisted generalized complex structure is a twisted generalized complex
submanifold.  We also prove that a twisted generalized complex manifold has a natural
Poisson structure.  We also discuss generalized K\"ahler submanifolds.
\footnote{MSC: 53C56, 53D17, 53D35} \footnote{Keywords: generalized complex geometry, Poisson bivector, Poisson-Dirac submanifold}
\end{abstract}
%=========================================================================================
%               START OF DOCUMENT.
%=========================================================================================

%%=========================================================================================
%\section*{Abstract}
%%=========================================================================================
%We introduce the notion of twisted generalized complex submanifolds and describe an
%equivalent characterization in terms of Poisson-Dirac submanifolds.  Our characterization
%recovers a result of Vaisman \cite{VAIS}.  An equivalent characterization is also given
%in terms of spinors.  As a consequence, we show that the fixed locus of an involution
%preserving a twisted generalized complex structure is a twisted generalized complex
%submanifold.  We also prove that a twisted generalized complex manifold has a natural
%Poisson structure.  We also discuss generalized K\"ahler submanifolds.

%=========================================================================================
\section{Introduction}
%=========================================================================================

Throughout this paper \(M\) will denote a smooth manifold.
Generalized complex structures were originally defined by Hitchin \cite{HIT1}, and further
studied by Gualtieri \cite{GUALT_thesis}.  Examples of generalized complex structures
include symplectic and complex manifolds.  In order to define generalized complex
structures we will recall some structures on \(\gtb{M}\).  The {\it Courant bracket} was
defined in \cite{COUR} as
\begin{equation}\label{cb}
  \cb{X+\xi}{Y+\eta}=\lb{X}{Y}+\ld{X}\eta-\ld{Y}\xi-\tfrac{1}{2}d\big(\eta(X)
  -\xi(Y)\big)
\end{equation}
for all \(X,Y\in\ovf{M}\) and \(\xi,\eta\in\dof{M}\).  There also exist smoothly varying
nondegenerate symmetric bilinear forms on each fibre of \(\gtb{M}\).  These are defined
as
\begin{equation}\label{met}
\ip{X+\xi}{Y+\eta}=\tfrac{1}{2}\big(\xi(Y)+\eta(X)\big)
\end{equation}
for all \(X,Y\in\tsp{m}{M}\), \(\xi,\eta\in\cts{m}{M}\) and \(m\in M\).

A {\it generalized complex structure}
is a smooth map \fct{J}{\gtb{M}}{\gtb{M}} such that \(J^2=-\id\), \(JJ^*=\id\), and the
\peb of \(J\) is involutive with respect to the Courant bracket.

The primary objects of study in this paper are twisted manifolds.
A manifold \(M\) endowed with a
closed 3-form \(\Omega\) will be called {\it twisted}.  A twisted manifold has another
well known bracket on \(\gvf{M}\): the {\it twisted Courant bracket}.
This bracket was defined in \cite{SEV-WEIN} as
\begin{equation}\label{cb-tw}
  \twbr{X+\xi}{Y+\eta}=\lb{X}{Y}+\ld{X}\eta-\ld{Y}\xi-\tfrac{1}{2}d\big(\eta(X)
  -\xi(Y)\big)+\is{Y}\is{X}\Omega.
\end{equation}

A {\it twisted generalized complex structure} is a smooth map \fct{J}{\gtb{M}}{\gtb{M}}
such that \(J^2=-\id\), \(JJ^*=\id\) and the \peb of \(J\) is involutive with respect to
\eqref{cb-tw}, rather than \eqref{cb}.
The triple \((M,\Omega,J)\) will be called a {\it twisted generalized complex manifold}.
Indeed, generalized complex manifolds can be defined for arbitrary Courant algebroids.
See \cite{BART_thesis} for details.

The aim of this work is to characterize when a submanifold of a twisted generalized
complex manifold is also a twisted generalized complex manifold.  In the untwisted case,
several notions of generalized complex submanifolds have been recently introduced.  The
notion defined here is similar to the one in \cite{BEN-BOY} and \cite{VAIS}.  A
different notion of generalized complex submanifolds appears in \cite{GUALT_thesis}.

\begin{defn}
A {\it twisted immersion}, from one twisted manifold \((N,\Upsilon)\) to another
\((M,\Omega)\), is defined as a smooth immersion \fct{h}{N}{M} with
\(\Upsilon=h^*\Omega\).  A {\it twisted generalized complex immersion} from
\((N,\Upsilon,J')\) to \((M,\Omega,J)\) is a twisted immersion 
\fct{h}{(N,\Upsilon)}{(M,\Omega)} such that the pullback of the \peb of \(J\) is the
\peb of \(J'\).  In this case, \(N\) is called a {\it twisted generalized complex
submanifold} of \(M\).
\end{defn}

By splitting vectors and covectors, a twisted generalized complex
structure can be written as
\begin{equation}\label{split} J=\bms{\phi&\pis\\\sigf&-\phi^*}.\end{equation}
Here \(\phi\) is an endomorphism of \(\tb{M}\),
\fct{\pis}{\ctb{M}}{\tb{M}} is the bundle map induced by a bivector field \(\pi\), and
\fct{\sigf}{\tb{M}}{\ctb{M}} is the bundle map induced by a 2-form \(\sigma\).
The fact that \(J^2=-\id\) also leads to the following formulas:
\begin{equation}\label{square}
\phi^2+\pis\sigf=-\id, \quad \phi\pis=\pis\phi^*, \quad\text{and}\quad
    \phi^*\sigf=\sigf\phi.
\end{equation}
These facts, and others, were first noted in \cite{CRAIN}.  These results were also
described using Poisson quasi-Nijenhuis manifolds in \cite{STIEN-XU2}.
The conditions for a submanifold to be twisted generalized complex will
be expressed in terms of this splitting.
Using the theory of Lie bialgebroids we also show that \(\pi\) from \eqref{split}
is a Poisson bivector field.  Which, for the untwisted case, is a standard result
\cite{AB-BOY,CRAIN}.
%%%  BAS-BOY has nothing, GUALT_thesis also doesn't.

Our work was inspired by \cite{STIEN-XU}, where
the reduction of generalized complex structures is studied. The main result was also
independently obtained by Vaisman \cite{VAIS}.

This paper is organized as follows.  In Section 2, we recall some of the basic facts of
Dirac structures.  In particular we describe the pull back.
In Section 3, we prove that a twisted generalized complex manifold carries a
natural Poisson structure.  In Section 4, we define the induced generalized complex
structure, and characterize when it has the required properties.  In Section 5, we
prove the main theorem of this paper, and provide examples.  Twisted
generalized complex involutions are also introduced in this section.  
In Section 6, we determine when a submanifold of a holomorphic Poisson manifold 
is itself endowed with an induced holomorphic Poisson structure.  Section 7 is a 
restatement of our main result in terms of spinors.  The last section discusses
generalized K\"ahler submanifolds.

\subsection*{Acknowledgements}
%The authors would like to thank Ping Xu for useful discussions.  
We would 
%also 
like to thank Izu Vaisman for pointing out an error in an earlier version of this paper.

%=========================================================================================
\section{Dirac structures}\label{ds}
%=========================================================================================
The aim of this section is to recall Dirac structures, and their
pull backs.  Before considering bundles, we will consider
a vector space \(V\).  In this case a {\it Dirac structure} is nothing more than a
maximal isotropic subspace of \(\dvs\).  Let \(q_1\) denote the projection of \(\dvs\)
onto \(V\), and \(q_2\) the projection onto \(V^*\).

If \(L\) is a Dirac structure then there exists a natural skew-symmetric bilinear form
\(\Lambda\) on \(L\) defined by
\[\Lambda(X+\xi,Y+\eta)=\xi(Y)=-\eta(X)\text{ for all }X+\xi,Y+\eta\in L.\]
It is easy to see that
\begin{align*}
  \Lambda(X+\xi_1,Y+\eta_1)
    &=\Lambda(X+\xi_2,Y+\eta_2)\text{ for all }X+\xi_{1,2},Y+\eta_{1,2}\in L,
  \intertext{and}
  \Lambda(X_1+\xi,Y_1+\eta)
    &=\Lambda(X_2+\xi,Y_2+\eta)\text{ for all }X_{1,2}+\xi,Y_{1,2}+\eta\in L.
\end{align*}
Hence, there exists a 2-form \(\ep\) on \(q_1(L)\) defined by
\[\ep(X,Y)=\Lambda(X+\xi,Y+\eta)\text{ for all }X+\xi,Y+\eta\in L,\]
and a 2-form \(\theta\) on \(q_2(L)\) defined by
\[\theta(\xi,\eta)=-\Lambda(X+\xi,Y+\eta)\text{ for all }X+\xi,Y+\eta\in L.\]

If $X\in q_1(L)$ then there exists some $\xi\in V^*$ with $X+\xi\in L$; furthermore
\(\ep(X,Y)=\xi(Y)\) for all \(Y\in q_1(L)\).  Thus \(\is{X}\ep=\xi\res{q_1(L)}\),
and
\[X+\xi\in L \quad\iff\quad X\in q_1(L)\text{ and }\is{X}\ep=\xi\res{q_1 (L)}.\]
Thus knowing the Dirac structure $L$ is exactly the same as knowing the subspace
\(q_1(L)\) and the 2-form \(\ep\).
Similarly, $L$ is equivalent to the pair $(q_2(L), \theta)$.  Thus any subspace
$R\se V$ endowed with a 2-form $\ep$ on \(R\) defines a Dirac structure \(L(R,\ep)\):
\[L(R,\ep)=\set{X+\xi\in R\ds V^*:\is{X}\ep=\xi\res{R}},\]
and any subspace $S\se V^*$ endowed with a 2-form $\pi$ on \(S\) defines a Dirac
structure \(L(S,\theta)\):
\[L(S,\theta)=\set{X+\xi\in V\ds S:\theta(\xi,\eta)=-\eta(X)\text{ for all }\eta\in S}.\]

Details of these constructions can be found in \cite{COUR}.  Let \(W\) be another
vector space and \fct{\varphi}{V}{W} a linear map.  The map \(\varphi\) can be used to
both  pull Dirac structures back from \(W\) to \(V\), and push Dirac structures forward 
from \(V\) to \(W\).  Let \((R,\ep)\) be a Dirac structure on \(W\), with \(R\se W\) and
\(\ep\in\kfo{2}{R}\).  A Dirac structure on \(V\) is defined by
\((\varphi^{-1} R,\varphi^*\ep)\).  This Dirac structure is called the {\it pull back}
of \((R,\ep)\) under $\varphi$.  Similarly if \((S,\theta)\) is a Dirac structure on
\(V\), with \(S\se V^*\) and \(\theta\) defined on \(S\), then 
\(((\varphi^*)^{-1} S,\varphi_*\theta)\) defines a Dirac structure on \(W\).  This Dirac
structure is
called the {\it push forward} of \((S,\theta)\) under $\varphi$.  These two maps of Dirac
structures are denoted by \(\F_\varphi\) and \(\B_\varphi\).  It is very easy to see that
for a Dirac structure \(L\) on \(W\)
\begin{align*}
B_\varphi(L)=\set{X+\varphi^*\xi:X\in V,\xi\in W^*\text{ such that }\varphi X+\xi\in L'},
\intertext{and for a Dirac structure \(L'\) on \(V\)}
\F_\varphi(L')=\set{\varphi X+\xi:X\in V,\xi\in W^*\text{ such that }X+\varphi^*\xi\in L}.
\end{align*}

Dirac structures can also be defined for a twisted manifold \((M,\Omega)\).  A
\emph{Dirac structure} is a smooth subbundle \(L\se\gtb{M}\) for which each fibre is a
Dirac structure of the corresponding fibre of $TM\oplus T^*M$, and whose space of
sections is closed under the twisted Courant bracket \eqref{cb-tw}.  The restriction of
the twisted Courant bracket to a Dirac structure is a Lie bracket; thus a Dirac structure
is naturally a Lie algebroid.

The definitions of push forward and pull back can be reformulated for Dirac structures on
manifolds.  We will only consider the pull back of a Dirac structure, but more on
the push forward can be found in \cite{BUR-RAD} and \cite{STIEN-XU}.  We note that the
pull back of a Dirac structure is automatically a maximal isotropic, but it need not be
smooth or involutive.

The last lemma of this section will be used to characterize when the
pullback bundle is involutive.  Let \((M,\Omega)\) and \((N,\Upsilon)\) be two twisted
manifolds with an immersion \fct{\varphi}{N}{M}.  Two sections
\(\sig_N=X+\xi\in\gvf{N}\) and \(\sig_M=Y+\eta\in\gvf{M}\)  are said to be
\emph{\(\varphi\)-related}, denoted by \(\frel{\varphi}{\sig_N}{\sig_M}\), if
\(Y=\varphi_*X\) and \(\xi=\varphi^*\eta\).  The following lemma is an extension of
Lemma 2.2 from \cite{STIEN-XU}.  This lemma is also true for complex sections, which
is when it will be applied.

\begin{lem}\label{15}
Assume that \(\sig_N^i\in\osn{\gtb{N}}\) and \(\sig_M^i\in\osn{\gtb{M}}\) satisfy
$\frel{\varphi}{\sig_N^i}{\sig_M^i}$, for $i=1,2$.  Then, if \(\varphi\) is a twisted
immersion,
\[\frel{\varphi}{\cb{\sig_N^1}{\sig_N^2}_{\Upsilon}}{\twbr{\sig_M^1}{\sig_M^2}}.\]
\end{lem}

\begin{proof}
Write $\sig_N^i=X^i+\xi^i$ and $\sig_M^i=Y^i+\eta^i$, where $X^i+\xi^i\in\gvf{N}$
and $Y^i+\eta^i\in\gvf{M}$, for $i=1,2$. Since $\frel{\varphi}{\sig_N^i}{\sig_M^i}$,
for \(i=1,2\), then $\varphi_*X^i=Y^i$ and $\varphi^*\eta^i=\xi^i$.  By definition
\[ \cb{\sig_N^1}{\sig_N^2}_{\Upsilon}=\lb{X^1}{X^2}+\ld{X^1}\xi^2-\ld{X^2}\xi^1
    +\tfrac{1}{2}d\big(\xi^1(X^2)-\xi^2(X^1)\big)+\is{X^2}\is{X^1}\Upsilon, \]
and
\[ \twbr{\sig_M^1}{\sig_M^2}=\lb{Y^1}{Y^2}+\ld{Y^1}\eta^2-\ld{Y^2}\eta^1
    +\tfrac{1}{2}d\big(\eta^1(Y^2)-\eta^2(Y^1) \big)+\is{Y^2}\is{Y^1}\Omega .\]
Now
\[ \varphi_*\lb{X^1}{X^2}=\lb{\varphi_*X^1}{\varphi_*X^2}=\lb{Y^1}{Y^2} ,\]
and
\begin{align*}
\varphi^*(\ld{Y^1}{\eta^2}) = & \varphi^*(i_{Y^1}d\eta^2+di_{Y^1}\eta^2)  
  = \varphi^* i_{\varphi_* X^1}d\eta^2+\varphi^* d\big(\eta^2(Y^1)\big) \\ 
  =& i_{X^1}\varphi^* d\eta^2+d\big(\xi^2(X^1)\big)
  = i_{X^1}d\xi^2+d i_{X^1}\xi^2
  = \ld{X^1}\xi^2 .\end{align*}
Similarly,
$$\varphi^*(\ld{Y^2}{\eta^1}) = \ld{X^2}\xi^1 .$$
The second last term becomes
$$\varphi^*d\big(\eta^1(Y^2)-\eta^2(Y^1)\big)=d\big(\xi^1(X^2)-\xi^2(X^1)\big),$$
since $$\varphi^*\big(\eta^1(Y^2)\big)=\varphi^*\big(\eta^1(\varphi_*X^2)\big)=
(\varphi^*\eta^1)(X^2)=\xi^1(X^2) .$$
Finally, because \(\varphi\) is a twisted immersion, the following holds.
\[\varphi^*\is{Y^2}\is{Y^1}\Upsilon=\varphi^*\is{\varphi_*X^2} \is{\varphi_*X^1}\Upsilon
  =\is{X^2}\is{X^1}\varphi^*\Upsilon=\is{X^2}\is{X^1}\Omega.\]
\end{proof}

%=========================================================================================
\section{The Poisson bivector field associated to a generalized complex structure} 
%=========================================================================================

For the usual Courant bracket it is well known that the existence of a generalized
complex structure leads to a Poisson bivector \cite{AB-BOY,CRAIN, ZAB-ETC, HU}.  In this
chapter we will obtain the same result for arbitrary Courant algebroids and also give a
new way of expressing the Poisson bivector.

A \emph{Courant algebroid} \cite{LIU-WEIN-XU} is a triple consisting of a vector bundle
$E\to M$ equipped with a non degenerate symmetric bilinear form $\ip{\cdot}{\cdot}$, a
skew-symmetric bracket $\cb{\cdot}{\cdot}$ on $\osn{E}$, and a smooth bundle map
$E\xrightarrow{\rho}TM$ called the anchor.  These induce a natural differential operator
$\DD :\sfn{M}\to\osn{E}$ defined by
\begin{equation} \label{c1} \ip{\DD f}{a}=\thalf\rho(a)f ,\end{equation}
for all $f \in \sfn{M}$ and $a\in\osn{E}$.  These structures must obey the following
formulas for all \(a,b,c\in\osn{E}\) and \(f,g\in \sfn{M}\),
\begin{gather}
  \anchor(\cb{a}{b})=\lb{\anchor(a)}{\anchor(b)}, \label{c2} \\
  \lcb{\cb{a}{b}}{c}+\lcb{\cb{b}{c}}{a}+\lcb{\cb{c}{a}}{b}=\tfrac{1}{3}
    \DD\big(\ip{\cb{a}{b}}{c}+\ip{\cb{b}{c}}{a}+\ip{\cb{c}{a}}{b}\big), \label{c3} \\
  \cb{a}{f b}=f\cb{a}{b}+\big(\anchor(a)f\big)b-\ip{a}{b}\DD f, \label{c4} \\
  \anchor\circ\DD=0,\text{ i.e. }\ip{\DD f}{\DD g}=0, \label{c5} \\
  \anchor(a)\ip{b}{c}
    =\ip{\cb{a}{b}+\DD \ip{a}{b}}{c}+\ip{b}{\cb{a}{c}+\DD \ip{a}{c}} \label{c6}.
\end{gather}

The relation \begin{equation} \label{c7} \cb{\DD f}{a}+\DD\ip{\DD f}{a}=0 \end{equation}
is a consequence of these conditions \cite{ROYT_thesis}.  

A smooth subbundle $L$ of a Courant algebroid is called a \emph{Dirac subbundle} if it is
a maximal isotropic, with respect to $\eip{\cdot}{\cdot}$, and its space of sections
$\osn{L}$ is closed under $\cb{\cdot}{\cdot}$.  While not all Courant algebroids are Lie
algebroids (since the Jacobi identity is not satisfied), their Dirac subbundles are Lie
algebroids.

\begin{eg}[\cite{COUR}] \label{std}
Given a smooth manifold $M$, the bundle $\tb{M}\to M$ carries a natural Courant algebroid
structure, where the anchor is the identity map and the pairing
and bracket are given, respectively, by
\begin{gather*}
  \ip{X+\xi}{Y+\eta}=\tfrac{1}{2}\big(\xi(Y)+\eta(X)\big), \\
  \cb{X+\xi}{Y+\eta}
    =\lb{X}{Y}+\ld{X}\eta-\ld{Y}\xi+\tfrac{1}{2}d\big(\xi(Y)-\eta(X)\big), 
\end{gather*}
$\forall X,Y \in \vf{M}$, $\forall \xi,\eta \in \dkf{1}{M}$.
\end{eg}

Let $E$ be a Courant algebroid on a smooth manifold $M$.  And let
\[ \xymatrix{
E \ar[d] \ar[r]^J & E \ar[d] \\ 
M \ar[r]_{\id} & M 
} \]
be a vector bundle map such that $J^2=-\id$.  Then the complexification
$E_\C:=E\otimes\C$ --- with the extended $\C$-linear Courant algebroid structure ---
decomposes as the direct sum $L\oplus\cjg{L}$ of the eigenbundles of $J$.  Here $L$ is
associated to the eigenvalue $+i$ and its complex conjugate $\cjg{L}$ to $-i$.  The
bundle map $J$ is called a \emph{generalized complex structure} if $J$ is orthogonal
with respect to $\ip{\cdot}{\cdot}$ --- this forces $L$ and $\cjg{L}$ to be isotropic
--- and the spaces of sections $\osn{L}$ and $\osn{\cjg{L}}$ are closed under the
Courant bracket, or equivalently, $J$ is ``integrable'':
\[ \cb{Jx}{Jy}-\cb{x}{y}-J\big(\cb{Jx}{y}+\cb{x}{Jy}\big)=0, \qquad \forall x,y\in\osn{E} .\]

Since the pairing is non degenerate, the map
\[ \Xi:E\xrightarrow{\isom}E^*:e\mapsto\ip{e}{\cdot} \]
is an isomorphism of vector bundles.
One has: $\Xi^*=\Xi$ (modulo the canonical isomorphism $(E^*)^*=E$) and
$\Xi\circ J+J^*\circ\Xi=0$.

\begin{prop}
The bracket
\begin{equation}
  \label{c8} \pb{f}{g}=2\ip{J\DD f}{\DD g}, \qquad f,g\in\sfn{M}
\end{equation}
is a Poisson structure on the manifold $M$.
\end{prop}

\begin{proof}
It is easy to see that this bracket is a skew-symmetric derivation of $\sfn{M}$. It
remains to check that the Jacobi identity is satisfied.  Since $J$ is integrable, we have
\[ \cb{J\DD f}{J\DD g}-\cb{\DD f}{\DD g}-J(\cb{J\DD f}{\DD g}+\cb{\DD f}{J\DD g})=0 ,\]
for all $f,g\in\sfn{M}$.  Pairing with $\DD h$, we obtain
\begin{equation} \label{c9} \ip{\cb{J\DD f}{J\DD g}}{\DD h}-\ip{\cb{\DD f}{\DD g}}{\DD h}
  -\ip{J(\cb{J\DD f}{\DD g}+\cb{\DD f}{J\DD g})}{\DD h}=0.
\end{equation}
We compute the first term of \eqref{c9}:
\begin{align*}
& \ip{\cb{J\DD f}{J\DD g}}{\DD h} && \\ 
=& \ip{\cb{J\DD f}{J\DD g}+\DD\ip{J\DD f}{J\DD g}}{\DD h} 
&& \text{by \eqref{c5}} \\ 
=& \anchor(J\DD f)\ip{J\DD g}{\DD h}- 
\ip{J\DD g}{\cb{J\DD f}{\DD h}+\DD\ip{J\DD f}{\DD h}} 
&& \text{by \eqref{c6}} \\ 
=& \;2\ip{\DD\ip{J\DD g}{\DD h}}{J\DD f}
-\ip{J\DD g}{\DD\ip{\DD h}{J\DD f}+\DD\ip{J\DD f}{\DD h}} 
&& \text{by \eqref{c1} and \eqref{c7}} \\ 
=& \;2\pb{f}{\pb{g}{h}}-2\pb{g}{\pb{f}{h}} 
&& \text{by \eqref{c8}} \\ 
=& \;2\pb{f}{\pb{g}{h}}+2\pb{g}{\pb{h}{f}} 
&& 
\end{align*}
By \eqref{c1}, \eqref{c2} and \eqref{c5}, the second term of \eqref{c9} vanishes.
Finally, the third term of \eqref{c9} gives:
\begin{align*} 
&-\ip{J(\cb{J\DD f}{\DD g}+\cb{\DD f}{J\DD g})}{\DD h} && \\ 
=&-\ip{J(\DD\ip{\DD g}{J\DD f}-\DD\ip{\DD f}{J\DD g})}{\DD h} 
&& \text{by \eqref{c7}} \\
=& -\pb{\pb{f}{g}}{h}+\pb{\pb{g}{f}}{h} 
&& \text{by \eqref{c8}} \\ 
=& \;2\pb{h}{\pb{f}{g}} && 
\end{align*}
\end{proof}

\begin{prop}
Let $\Pi$ denote the bivector field on $M$ associated to the Poisson bracket
\eqref{c8}. We have
\[ \Pi\diese=\thalf\anchor\circ J\circ\Xi^{-1}\of\anchor^*, \]
where $\Pi\diese:T^*M\to TM$ is the vector bundle map equivalent to
$\Pi\in\osn{\wedge^2 T^*M}$ via
$\Pi(\alpha,\beta)=\beta(\Pi\diese\alpha),\;\forall\alpha,\beta\in\dkf{1}{M}$.
The Hamiltonian vector field associated to the function $f\in\sfn{M}$ is
\[ X_f=\Pi\diese df=\anchor\of J\of\DD f .\]
And the characteristic distribution is
\[ \Pi\diese(\ctb{M})=\anchor J (\ker\anchor)^\perp, \]
where $(\ker\anchor)^\perp$ refers to the subbundle of $E$ orthogonal to $\ker\anchor$
with respect to $\ip{\cdot}{\cdot}$.
\end{prop}

\begin{proof}
One has
\[ X_f(g)=\pb{f}{g}=2\ip{J\DD f}{\DD g}=\anchor(J\DD f)\big(g\big) \]
and
\[ \Pi\diese(df)=X_f=\anchor J\DD f=\thalf\anchor J\Xi^{-1}\anchor^*(df) \]
since \eqref{c1} can be reinterpreted as
$\ip{\DD f}{\cdot}=\thalf\Xi^{-1}\anchor^* df$.
\end{proof}

\begin{prop}
If $E=\tb{M}$ is the standard Courant algebroid of Example~\ref{std}
and the matrix representation of $J$ relative to the above direct sum decomposition is
$\big(\begin{smallmatrix} \phi & \pi\diese \\ \sigf & -\phi^* 
\end{smallmatrix}\big)$, then $\Pi=\pi$.
\end{prop}

\begin{proof}
Here $\DD$ coincides with the de Rham differential $d$. Thus
\[ \pb{f}{g}_{\Pi}=2\ip{J\DD f}{\DD g}
=2\ip{\pi\diese df - \phi^*df}{dg}=dg(\pi\diese df)=\pb{f}{g}_{\pi} .\]
\end{proof}

Recall that the complexification of $E$ decomposes as the direct sum $E_\C=L_+\oplus L_-$, where $L_{\pm}$ are Dirac structures (with anchor maps $\anchor_{\pm}$). Thus $(L_+,L_-)$ is a complex Lie bialgebroid \cite{MAC-XU}, where $L_{\pm}^*$ is identified with $L_{\mp}$ via $\Xi$.  As shown in \cite[Proposition 3.6]{MAC-XU}, to any complex Lie bialgebroid is associated a complex bivector field $\varpi$ on $M$ given by
\[ i\varpi^\sharp=\anchor_-\of\Xi^{-1}\of\anchor_+^*
   = -\anchor_+\of\Xi^{-1}\of\anchor_-^* .\]

\begin{lem}
The Poisson bivector $\varpi$ coming from the Lie bialgebroid structure $(L_+,L_-)$ is real and coincides with $\Pi$.
\end{lem}

\begin{proof}
It suffices to observe that the following two compositions are both equal to $\varpi\diese$:
\begin{gather*}
T_\C M\xleftarrow{\anchor_-}L_-\xleftarrow{-i\cdot}L_-\xleftarrow{\Xi^{-1}} 
L_+^*\xleftarrow{\anchor_+^*}T_\C^* M \\ 
T_\C M\xleftarrow{\anchor_+}L_+\xleftarrow{+i\cdot}L_+\xleftarrow{\Xi^{-1}} 
L_-^*\xleftarrow{\anchor_-^*}T_\C^* M 
\end{gather*}
and that their sum
\[ T_\C M\xleftarrow{\anchor}E_\C\xleftarrow{J}E_\C\xleftarrow{\Xi^{-1}} 
E_\C^*\xleftarrow{\anchor^*}T_\C^* M \]
is equal to $2\Pi\diese$.
\end{proof}

\begin{prop}
$\anchor(L_+)\cap\anchor(L_-)=\Delta\otimes\C$ with $\Delta=\anchor(J\ker\anchor)$
\end{prop}

\begin{proof}
If $v\in\anchor(L_+)\cap\anchor(L_-)$, then
$\cjg{v}\in\anchor(L_-)\cap\anchor(L_+)$.
Hence there exists a subbundle $\Delta$ of $TM$ such that
$\anchor(L_+)\cap\anchor(L_-)=\Delta\otimes\C$.
For all $k\in\ker\anchor$, one has
$\anchor(L_+)\ni\anchor(\tfrac{1+iJ}{2}k)=\tfrac{i}{2}\anchor(Jk)$.
Therefore, $\anchor(J\ker\anchor)\subset\anchor(L_+)$.
Since $J$ is a real, $\anchor(J\ker\anchor)\subset\anchor(L_+)\cap\anchor(E)=\Delta$.
It remains to prove the converse inclusion: $\Delta\subset\anchor(J\ker\anchor)$.
Since $\Delta=\anchor(L_+)\cap\anchor(E)$, given $\delta\in\Delta$, there exists $l_+\in L_+$ such that $\anchor(l_+)=\delta=\anchor(\cjg{l_+})$.
Thus $\delta=\rho(\tfrac{l_+ +\cjg{l_+}}{2})= 
\rho\big(J(\tfrac{l_+ -\cjg{l_+}}{2i})\big)$ with $\tfrac{l_+ -\cjg{l_+}}{2i}\in\ker\anchor$.
\end{proof}

\begin{rem}
If $E=\tb{M}$ is the standard Courant algebroid of Example~\ref{std}, then
$(\ker\anchor)^\perp=T^*M=\ker\anchor$.
Therefore, in this particular case, $\Pi\diese(T_\C^*M)=
\anchor(L_+)\cap\anchor(L_-)$, recovering Gualtieri's result \cite{GUALT_thesis}.

It would be interesting to explore when the symplectic foliation $\Pi\diese(T_{\mathbb{C}}^*M)$ coincides with 
$\anchor(L_+)\cap\anchor(L_-)$ for arbitrary Courant algebroids.
\end{rem}

%=========================================================================================
\section{The induced generalized complex structure}\label{ind}
%=========================================================================================

Consider two twisted manifolds \((M,\Omega)\) and \((N,\Upsilon)\) with an immersion
\fct{h}{N}{M}.  Also, assume that there is a generalized complex structure \(J\) on
\(M\) with eigenbundles \(L_+\) and \(L_-\).  The goal of this section is to characterize
when the pull backs of \(L_+\) and \(L_-\) give a generalized complex structure on \(N\).
The pull backs of \(L_+\) and \(L_-\) will be called the {\it induced bundles}, and
are given by
\[L'_\pm=\B_h(L_\pm)=\set{X+h^*\xi:X\in\tbc{N},\xi\in\ctbc{M}
  \text{ such that }h_*X+\xi\in L_\pm},\]
By definition, both \(L_+'\) and \(L_-'\) are maximal isotropics, but they need not be
smooth or involutive subbundles.  The bundles may also have nontrivial intersection.  The
rest of this section is devoted to characterizing when the induced bundles have the
desired properties.  The first of these properties to be addressed will be the
intersection property.

Because \(L_+'\) and \(L_-'\) are both maximal isotropics,
it suffices to check that they span \(\gtbc{N}\).  Consider the subbundle
\(B=\tb{N}\ds\ctb{M}\res{N}\) of \(\gtb{M}\).  Its orthogonal, \(B^\perp=\tb{N}^o\), is
the kernel of the natural projection \fct{s}{B}{\gtb{N}}, which maps
\mapto{X+\xi}{X+h^*\xi}.  It is not hard to see that
\(s((B\cap JB)\otimes\C)=L'_++L'_-\).  Thus the decomposition,
\(\gtbc{N}= L'_+ \ds L'_-\), holds if and only if \(B=B\cap JB + B^\perp\).  The
preceding can be summarized in the following proposition.

\begin{prop}\label{1}
The following assertions are equivalent.
\begin{enumerate}
  \item The subbundle $L'_+$ is the $+i$-eigenbundle of a -- not necessarily smooth --
    automorphism $J'$ of $\gtb{N}$ such that ${J'}^2=-\id$ and $J'J'^*=\id$.
  \item  \(B=B\cap JB + B^\perp\).
  \item  \(JB\se B+JB^\perp\).
  \item  \(JB^\perp\cap B\se B^\perp\).
\end{enumerate}
\end{prop}

Conditions (3) and (4) follow from elementary calculations.  In the sequel we will
assume that the assertions of \rprop{1}
are satisfied.  Consider the restriction of \(J\) and \(s\) to the $J$-invariant subspace
$B\cap JB$; the latter map will be denoted by \(s'\).  The kernel of \(s'\)  is
\(B^\perp\cap JB\).  Under \(J\), this kernel is mapped to \(JB^\perp\cap B\).  This
must be in \(JB\cap B\) and also, by \rprop{1}, in \(B^\perp\), however \(B^\perp\se B\).
So the image of the kernel of \(s'\) is in \(B^\perp\cap JB\cap B=B^\perp\cap JB\).  Thus
the kernel of \(s'\) is \(J\)-invariant and $J|_{B\cap JB}$ induces an automorphism of
\(\gtb{N}\):
\begin{equation}\label{11} \xymatrix{ B\cap JB \ar[d]_s \ar[rr]^J && B\cap JB \ar[d]^s \\
  \gtb{N} \ar[rr]_{J'} && \gtb{N}.} \end{equation}
The induced automorphism is nothing but $J'$ from Proposition~\ref{1}. Indeed, the
complexification of the above commutative diagram gives
\[\xymatrix{ (L_+ \cap B_\C)\ds(L_-\cap B_\C) \ar[d] \ar[rrr]^{(+i)\id\ds
  (-i)\id} &&& (L_+ \cap B_\C)\ds(L_-\cap B_\C) \ar[d] \\
  L'_+\ds L'_- \ar[rrr]_{(+i)\id\ds (-i)\id} &&& L'_+\ds L'_-.} \]

The following lemma relates condition (4) of \rprop{1} to the splitting of \(J\):
\[ J=\bms{\phi&\pis\\\sigf&-\phi^*}.\]

\begin{lem} \label{2}
The following assertions are equivalent.
\begin{enumerate}
\item \(JB^\perp\cap B\se B^\perp\).
\item \(\tb{N}\cap\pis(\tb{N}^o)=0\) and \(\phi(\tb{N})\se\tb{N}+\pis(\tb{N}^o)\).
\end{enumerate}
\end{lem}
\begin{proof}
The inclusion
\[ J(\tb{N}^o)\cap(\tb{N}\ds\ctb{M}\res{N})\se \tb{N}^o \] 
is true if, and only if, 
\[ \left\{ \begin{array}{l} \xi\in\tb{N}^o \\ 
J\xi\in\tb{N}\ds \ctb{M}\res{N} \end{array} \right\} 
\implies J\xi\in\tb{N}^o \]
if, and only if, 
\[ \left\{ \begin{array}{l} \xi\in\tb{N}^o \\ 
\pis\xi\in\tb{N} 
\end{array} \right\} 
\implies \left\{ \begin{array}{l} 
\pis\xi=0 \\ \phi^*\xi\in\tb{N}^o
\end{array} \right\} \]
if, and only if,
\[ \xi\in\tb{N}^o\cap(\pis)^{-1}(\tb{N})  
\implies \left\{ \begin{array}{l} 
\pis\xi=0 \\ \xi\in(\phi(\tb{N}))^o
\end{array} \right\} \]
if, and only if, 
\[ \pis(\tb{N}^o)\cap\tb{N}=0 \qquad\text{and}\qquad 
\tb{N}^o\cap(\pis)^{-1}(\tb{N})\se(\phi(\tb{N}))^o.\]
Since $\pi$ is skew-symmetric \((\pis)^{-1}(\tb{N})=(\pis(\tb{N}^o))^o\), and
\begin{align*}
\tb{N}^o\cap(\pis(\tb{N}^o))^o &\se(\phi(\tb{N}))^o\\
\tb{N}+\pis(\tb{N}^o) &\sse\phi(\tb{N}).
\end{align*}
\end{proof}

According to this lemma, the sum $\tb{N}+\pis (\tb{N}^o)$ must be direct.  In the sequel
\(\pr\) will denote the projection \fctto{\tb{N}\ds\pis(\tb{N}^o)}{\tb{N}}.  If \(\pi\)
is degenerate then neither the bundle \(\tb{N}\ds\pis(\tb{N}^o)\), nor the map
\(\pr\) are necessarily smooth.

For any \(\xi\in\ctb{N}\) we claim that if \(\eta,\eta'\in B\cap JB\) such that
\(\xi=h^*\eta=h^*\eta'\) then \(\pis\eta=\pis\eta'\).  Because \(\eta\)
and \(\eta'\) are preimages of \(\xi\) they differ by some element of \(\tb{N}^o\), and
as \(B\cap JB\) is \(J\)-stable both \(\pis\eta\) and \(\pis\eta'\) are in \(\tb{N}\).
However \(\tb{N}\cap\pis(\tb{N}^o)=\z\), and the difference of the two preimages is zero.
Thus the assignment $\xi\mapsto\pis\eta$ defines a skew-symmetric vector bundle map from
\(\ctb{N}\) to \(\tb{N}\).  Its associated bivector field on $N$ will be denoted by
$\pi'$.

The following technical lemmas will be used to show when \(J'\) is smooth.

\begin{lem} \label{4}
Let $X\in\tb{N}$ and $\xi\in\tb{N}^o$.  If \(\phi X+\pis\xi\in\tb{N}\) then
\(\phi X+\pis\xi=(\pr\of\phi) X\).
\end{lem}
\begin{proof}
For any \(X\in\tb{N}\) the second assertion of Lemma~\ref{2} gives \(\phi X=Y+\pis\eta\),
where \(Y\in\tb{N}\) and \(\eta\in\tb{N}^o\).  By definition \(Y=(\pr\of\phi)X\), and
\(\phi X+\pis\xi=Y+\pis(\eta+\xi)\).  Both \(Y\) and \(\phi X+\pis\xi\) are elements of
\(\tb{N}\); thus \(\pis(\eta+\xi)\) is also an element of \(\tb{N}\).  But
\(\tb{N}\cap\pis(\tb{N}^o)=\z\), and \(\eta+\xi\in\tb{N}^o\).  Thus \(\pis(\eta+\xi)=0\).
\end{proof}

\begin{lem} \label{6}
Let \fct{p_1}{\gtb{N}}{\tb{N}} and \fct{p_2}{\gtb{N}}{\ctb{N}} be the projections.
If $X\in\tb{N}$, then 
\begin{align}
(p_1 J')X &= (\pr\of\phi)X =\phi X +\pis\zeta \label{7}\intertext{and}
(p_2 J')X &= h^* (\sigf X-\phi^*\zeta), \label{8}
\end{align}
where $\zeta$ is some element of \(\tb{N}^o\) such that $X+\zeta\in B\cap JB$. 

If \(\xi\in\ctb{N}\), then 
\begin{align}
(p_1 J')\xi&= \pips\xi = \pis\eta \label{9}\intertext{and}
(p_2 J')\xi&=-(h^*\phi^*)\eta\label{10},
\end{align}
where $\eta$ is some element of $\ctb{M}\res{N}\cap B\cap JB$ such that $h^*\eta=\xi$.
\end{lem}
\begin{proof}
Consider $X\in\tb{N}$. Since \(s\) is surjective 
there exists some $\zeta\in\tb{N}^o$ such that $X+\zeta\in B\cap JB$ and $s(X+\zeta)=X$. 
Now $J(X+\zeta)=(\phi X+\pis\zeta)+(\sigf X-\phi^*\zeta)\in B$.  Therefore
$\phi X+\pis\zeta\in \tb{N}$ and, by Lemma~\ref{4}, $\phi X+\pis\zeta=(\pr\of\phi)X$.
Both \eqref{7} and \eqref{8} follow from \eqref{11}.

Now take \(\xi\in\ctb{N}\). Again, since \(s\) is surjective
there exists some \(\eta\in\ctb{M}\) such that \(\eta\in B\cap JB\) and \(s(\eta)=\xi\).
Now \(J(\xi)=\pis\eta-\phi^*\eta=\pips\xi-\phi^*\eta\), which is in \(B\).
Both \eqref{9} and \eqref{10} follow from \eqref{11}.
\end{proof}

For the remainder of this section, if $L$ is a smooth vector bundle then $\nss(L)$ will
denote the space of
all -- not necessarily smooth -- sections of $L$, and $\sms(L)$ the subspace of smooth
sections.

\begin{lem} \label{5}
Let \(\xi\in\nss(\tb{N}^o)\).
If \(\pis\xi\in\sms(\tb{M}\res{N})\), then \((h^*\phi^*){\xi}\in\sms(\ctb{N})\).
\end{lem}
\begin{proof}
As noted previously, if \(X\in\sms(\tb{N})\) then \(\phi X=Y+\pis\eta\), where
\(Y\in\osn{\tb{N}}\) and \(\eta\in\sms(\tb{N}^o)\). Now  
\((\phi^*\xi)(X)=\xi(\phi X)=\xi(Y)+\xi(\pis\eta)=\xi(\pis\eta)=-\eta(\pis\xi)\).
This function and its restriction to \(\tb{N}\) are smooth because $\pis\xi$ is.
\end{proof}

\begin{lem} \label{12}
Assume \(\pr\of\phi\) is a smooth map and $\eta\in\nss(\ctb{M}\res{N})$. If
$h^*\eta\in\sms(\tb{N})$ and $\pis\eta\in\sms(\tb{N})$ then
$(h^*\phi^*)\eta\in\sms(\ctb{N})$.
\end{lem}
\begin{proof}
Once again, if \(Y\in\sms(\tb{N})\) then $\phi Y=(\pr\of\phi)Y+\pis\zeta$ for some
$\zeta\in\sms(\tb{N}^o)$.  Now
\[ \big((h^*\phi^*)\eta\big)(Y)=(h^*\eta)(\phi Y)=(h^*\eta)\big((\pr\of\phi)Y\big)
  +(h^*\eta)(\pis\zeta)=(h^*\eta\big)((\pr\of\phi)Y\big)-\zeta\big(\pis(h^*\eta)\big). \]
Thus $((h^*\phi^*)\eta)(Y)$ is a smooth function, and the lemma follows.
\end{proof}

We are now ready to give the conditions \(J'\) must satisfy in order to be smooth.

\begin{prop} \label{13}
The vector bundle automorphism $J'$ of $\gtb{N}$ is smooth if, and only if,
\fct{\pr\of\phi}{\tb{N}}{\tb{N}} is smooth and $\pi'$ is a smooth bivector field on $N$.
\end{prop}
\begin{proof}
%\fbox{$\Rightarrow$}
First assume that $J'$ is smooth.  Thus $(p_1 J')X\in\sms(\tb{N})$ for all
$X\in\sms(\tb{N})$.  It follows from \eqref{7} that $(\pr\of\phi)$ must be smooth.
Also $(p_1 J')\xi\in\sms(\ctb{N})$ for all $\xi\in\sms(\ctb{N})$, and \eqref{9} shows
that $\pips$ is smooth.

Now for the other implication.
%\fbox{$\Leftarrow$}
For every $X\in\sms(\ctb{N})$ there is some $\zeta\in\nss(\tb{N}^o)$ such that \eqref{7} and
\eqref{8} are satisfied.  As $J$ is smooth both $\sigf$ and $\phi$ are smooth.
The smoothness of $\pr\of\phi$ and \eqref{7} show that
$\pis\zeta\in\sms(\tb{M}\res{N})$.
Thus, according to Lemma~\ref{5}, $(h^*\phi^*)\zeta\in\sms(\ctb{N})$, and the right
hand sides of \eqref{7} and \eqref{8} are smooth.  Finally
\(J'X=(p_1J')X+(p_2J')X\in\sms(\gtb{N})\).

Now take $\xi\in\sms(\ctb{N})$.
There must exist $\eta\in\nss(\ctb{M}\res{N})$ such that \eqref{9} holds, \eqref{10}
holds, and $h^*\eta=\xi$.
The smoothness of $\pi'$ and \eqref{9} show that $\pis\eta\in\sms(\tb{N})$.  Now
Lemma~\ref{12} gives $(h^*\phi^*)\eta\in\sms(\ctb{N})$, and the right hand
sides of \eqref{9} and \eqref{10} are smooth.  Finally,
\(J'\xi=(p_1J')\xi+(p_2J')\xi\in\sms(\gtb{N})\).
\end{proof}

We finish this section by using \rlem{15} to show when \(J'\) is integrable.

\begin{prop} \label{14}
If $J'$ is smooth then it is integrable.
\end{prop}
\begin{proof}
First, observe that the vector bundles $L_{\pm}\cap B_\C = (I\mp iJ)B_\C$ are
smooth.  Since $J'$ is smooth, its eigenbundles $L'_{\pm}$ are also smooth. It is not
hard to check that any smooth section of $L'_+$ is $h$-related to a smooth section of
$L_+\cap B_\C$.

Hence for any \(\sigma'_1,\sigma'_2\in\sms(L'_+)\) there exists
\(\sigma_1,\sigma_2\in\sms(L_+\cap B_\C)\) 
such that \(\frel{h}{\sigma_1}{\sigma'_1}\) and \(\frel{h}{\sigma_2}{\sigma'_2}\). 
Since $L_+$ is integrable $\twbr{\sigma_1}{\sigma_2}\in\sms(L_+)$, and it follows from
\rlem{15} that
\(\frel{h}{\twbr{\sigma_1}{\sigma_2}}{\icb{\sigma'_1}{\sigma'_2}}{\Upsilon}\).  Thus
\(\icb{\sigma'_1}{\sigma'_2}{\Upsilon}\in\sms(L'_+)\), and $L'_+$ is involutive with
respect to the \(\Upsilon\)-twisted bracket.
\end{proof}

%=========================================================================================
\section{Main theorem}\label{mt}
%=========================================================================================

The following definition will be used to characterize when a twisted submanifold is also generalized complex; see \cite{CRAIN-FERN} for the motivation of this definition.

\begin{defn}
Let \((M,\pi)\) be a Poisson manifold.  A smooth submanifold \(N\) of \(M\) is a
{\it Poisson-Dirac submanifold} of \(M\) if \(\tb{N}\cap\pis\kp=\z\), and the induced
Poisson tensor \(\pi'\) on \(N\) is smooth.
\end{defn}

The next theorem is the main result of this paper.  The untwisted version of this
result was obtained independently, using a different method, by Vaisman \cite{VAIS}.

\begin{thm}\label{main}
Let \((M,\Omega,J)\) be a twisted generalized complex manifold with
\(J=\bms{\phi&\pis\\\sigf&-\phi^*}\).  A twisted submanifold \(N\) of \(M\) inherits a
twisted generalized complex structure $J'$, making it a twisted generalized complex
submanifold, if and only the following conditions hold:
\begin{enumerate}
  \item  \(N\) is a Poisson-Dirac submanifold of \((M,\pi)\),
  \item  \(\phi(\tb{N})\se\tb{N}+\pis(\tb{N}^o)\),
  \item  \fct{\pr\of\phi}{\tb{N}}{\tb{N}} is smooth.
\end{enumerate}
The generalized complex structure \(J'\) on \(N\) is given by
\[ J'=\bms{\phi'&(\pi')^\sharp\\\sigf'&-(\phi')^*}. \]
Here \(\phi'=\pr\of\phi\res{\tb{N}}\), \(\pi'\) is the induced Poisson tensor, and
\begin{equation}\label{sigp}\sigf'(X)=h^*(\sigf X-\phi^*\zeta),\end{equation}
where \(\zeta\in(\tb{N})^o\) such that \(X+\zeta\in B\cap JB\), as in \rlem{6}
\end{thm}

\begin{proof}
This theorem is the construction and confirmation of the properties of \(J'\).  \rprop{1}
combined with \rlem{2} shows that \(J'^2=-\id\) and \(J'^*J'=\id\).  The smoothness of
\(J'\) follows from \rprop{13}, and the integrability of its \(+i\)-eigenbundle follows
from \rprop{14}.  The form of the generalized complex structure follows from \rlem{6}.
\end{proof}

For the following examples let \(\Omega=0\).

\begin{eg}\label{egc}
Let \((M,j)\) be a complex manifold, and let \(N\) be a smooth submanifold of \(M\).
There is a
generalized complex structure on \(M\) given by \(\phi=j\), \(\sigma=0\) and \(\pi=0\).
Because the Poisson structure is zero, \(N\) is automatically a Poisson-Dirac
submanifold.  Condition (b) of \rthm{main} becomes \(j(\tb{N})\se\tb{N}\), which is
exactly the requirement for \(N\) to be an immersed complex submanifold of \(M\).
Now \(\pr\of j=j\res{\tb{N}}\), which is a smooth map.  Thus \(N\) is a generalized
complex submanifold if, and only if, it is an immersed complex submanifold.
\end{eg}

\begin{eg}\label{egs}
Let \((M,\omega)\) be a symplectic manifold and \(N\) a smooth submanifold of \(M\).  The
generalized complex structure on \(M\) arising from \(\omega\) is given by \(\phi=0\),
\(\sigf=\omega_{\flat}\) and \(\pis=-\omega_{\flat}^{-1}\).  
Because \(\phi=0\), conditions (b)
and (c) of \rthm{main} are automatically satisfied.  Now \(N\) will be a Poisson-Dirac
submanifold of \(M\) if, and only if, \(N\) is a symplectic submanifold of \(M\).  Thus
\(N\) is a
generalized complex submanifold of \(M\) if, and only if, it is a symplectic submanifold.
\end{eg}

The last result of this section is an application of \rthm{main} to the stable locus of a
twisted generalized complex involution.  This result is similar to one for Poisson
involutions \cite{FERN-VAN,XU}.  Let \((M,\Omega,J)\) be a twisted generalized
complex manifold.  A {\it twisted generalized complex involution} is a diffeomorphism
\fct{\Psi}{M}{M} such that \(\Psi^2=\id\) and 
\begin{equation}\label{cplx}\Psi_*^*\of J=J\of\Psi_*^*.\end{equation}
Here \(\Psi^*_*\) is the map from \(\gtb{M}\) to \(\gtb{M}\) defined by
\(\Psi_*^*(X+\xi)=\Psi_*X+\Psi^*\xi\).
\begin{cor}
Let \((M,\Omega,J)\) be a twisted generalized complex manifold and let \(\Psi\) be a
twisted generalized complex involution of \(J\).  The fixed locus, \(N\), of \(\Psi\) is
a twisted generalized complex submanifold of \(M\).
\end{cor}
\begin{proof}
Let \(\xi\) be an arbitrary element of \(\ctb{M}\).  Equation \eqref{cplx} implies that
\((\Psi_*\pis\Psi^*)\xi=\pis\xi\).  Hence \(\Psi_*\pi=\pi\), and \(\Psi_*\) is a
Poisson involution.  Because \(\Psi_*\) is a Poisson involution, Proposition 4.1 of
\cite{XU} implies that \(N\) is a Dirac submanifold.  Thus \(N\) is a Poisson-Dirac
submanifold, and condition (a) of \rthm{main} is satisfied.

Take \(X\in\tb{N}\).  Equation \eqref{cplx} implies that
\(\Psi_*(\phi X)+\Psi^*(\sigf X)=\phi X+\sigf X\).  The vector field component of this
equality proves that \(\phi(\tb{N})\se\tb{N}\), and condition (b) of \rthm{main} is
satisfied.  Thus \(\pr\of\phi=\phi\res{\tb{N}}\), which is a smooth map.  Hence condition
(c) of \rthm{main} is satisfied.
\end{proof}

%=========================================================================================
\section{Holomorphic Poisson submanifolds}\label{ph}
%=========================================================================================

Let \((M,j,\pi)\) be a Poisson Nijenhuis manifold such that \fct{j}{\tb{M}}{\tb{M}} is
an integrable almost complex structure.  Such a structure is equivalent to a holomorphic
Poisson structure  The holomorphic Poisson tensor is
given by \(\Pi=\pi_j+i\pi\), where \(\pi_j^\sharp=\pis\of j^*\).

A generalized complex structure on \(M\) is given by, \cite{CRAIN,STIEN-XU2},
\begin{equation}\label{egpneq} J=\bms{j&\pis\\0&-j^*}. \end{equation}
In general, if \(N\) is a generalized complex submanifold then the induced generalized
complex need not have \(\sig'=0\).

Recall that \(\tb{N}\cap\pis(\tb{N}^o)=\z\) and \(\phi(\tb{N})\se\tb{N}+\pis(\tb{N}^o)\).
Thus, we can define the composition
\begin{equation}\label{124}
  \tb{N}\xrightarrow{\phi}\tb{N}\ds\pis(\tb{N}^o)\xrightarrow{\pr_2}\pis(\tb{N}^o).
\end{equation}

\begin{prop}
Consider the generalized complex structure \eqref{egpneq}, and let \(N\) be a generalized
complex submanifold of \(M\).  Now, \(\sig'=0\) if, and only if,
\begin{equation}\label{123}
\phi(\tb{N})\se\pis\big((\pr_2\of\phi(\tb{N}))^o\big).
\end{equation}
\end{prop}
\begin{proof}
Take \(X\in\tb{N}\).  Then \(\bms{X\\0}\in\tb{N}\ds\ctb{M}\res{N}=B\).  Since
\(B=B\cap JB+B^\perp\) and \(B^\perp=\tb{N}^o\) there exists a \(\zeta\in\tb{N}^o\)
such that \(\bms{X\\\zeta}\in B\cap JB\).  But then \(J\bms{X\\\zeta}\in B\) too, and
so \(\phi X+\pis\zeta\in\tb{N}\).

In other words, given \(X\in\tb{N}\) there exists \(\zeta\in\tb{N}^o\) such that
\(\phi X+\pis\zeta\in\tb{N}\).  Since \(\phi(\tb{N})\se\tb{N}\ds\pis(\tb{N}^o)\), this
is equivalent to \(\pr_2\of\phi(X)=-\pis\zeta\).  Recall that
\(\sigf'(X)=h^*(\sigf(X)-\phi^*\zeta)\).

Now, assume \(\sig=0\).  Then \(\sig'=0\) if, and only if \(\sigf'(X)=0\) for all
\(X\in\tb{N}\).  From the previous discussion, this will be true if, and only if
\(h^*(\phi^*\zeta)=0\) for all \(\zeta\in\tb{N}^o\) such that
\(\pis\zeta\in\pr_2\of\phi(\tb{N})\).  Which is equivalent to 
\begin{equation}\label{125}
  \tb{N}^o\cap(\pis)^{-1}(\pr_2\of\phi(\tb{N}))\se(\phi^*)^{-1}(\tb{N}^o).
\end{equation}

Let \(A=\pr_2\of\phi(\tb{N})\).  Now,
\[ (\pis)^{-1}(A)=\big((\pis)^*\big)^{-1}(A)=\big(\pis(A^o)\big)^o, \]
and
\[ (\phi^*)^{-1}(\tb{N}^o)=\big(\phi(\tb{N})\big)^o.\]
Hence, \eqref{125} becomes
\[\big(\tb{N}+\pis(A^o)\big)^o\se\big(\phi(\tb{N})\big)^o,\]
which is equivalent to 
\begin{equation}\label{eg1} \phi(\tb{N})\se\tb{N}+\pis(A^o), \end{equation}

So a generalized complex submanifold \(N\), of a generalized complex
manifold \(M\) with generalized complex structure \eqref{egpneq}, will have a generalized
complex structure of the same form as \eqref{egpneq} if and only if
\(\phi(\tb{N})\se\tb{N}+\pis(A^o)\).  Now, consider the following series of equivalent
statements:
\begin{align*}
  A &\se\pis(\tb{N}^o) \\
  A^o &\sse \big(\pis(\tb{N}^o)\big)^o \\
  A^o &\sse (\pis)^{-1}(\tb{N}) \\
  \pis(A^o)&\sse\pis\big((\pis)^{-1}(\tb{N})\big)=\tb{N}.
\end{align*}
Thus \eqref{eg1} becomes
\begin{equation} \phi(\tb{N})\se\pis(A^o). \end{equation}
\end{proof}

If \(N\) is both a complex submanifold of \((M,j)\), and a Poisson-Dirac submanifold of
\((M,\pi)\) then this condition will automatically be satisfied and there will be a
generalized complex structure of the form \eqref{egpneq} on \(N\).  Thus \(N\) will also
be a holomorphic Poisson manifold.

%=========================================================================================
\section{Spinors and generalized complex submanifolds}
%=========================================================================================

Generalized complex structures may also be realized using Clifford algebras and
spinors.  The aim of this section is to prove that generalized complex submanifolds can
also be realized using spinors.  Details for the definitions in this section can be
found in \cite{GUALT_thesis}, and
the sections of \cite{CHEV} cited therein.

Let \(V\) be a finite dimensional vector space and let \(\cliff\) denote the Clifford
algebra of \(\cdvs\).  There is an action of \(\cliff\) on \(\clf{V}\) defined by
\[ (X+\xi)\cdot\mu=\is{X}\mu+\xi\wedge\mu,\]
for all \(X+\xi\in\cdvs\) and \(\mu\in\clf{V}\).  The elements of \(\clf{V}\) are called
{\it spinors}.  Each spinor \(\mu\) has a {\it null space}:
\[ L_\mu=\set{X+\xi\in\cdvs:(X+\xi)\cdot\mu=0}. \]
This subspace is isotropic, and if it is also maximal isotropic then the spinor is
called {\it pure}.
Using the notation of \rsect{ds}, every pure spinor can be written as
\begin{equation}\label{spinorform} \mu=c(\det(R^o))\wedge\exp(\ep), \end{equation}
where \(c\in\C\) is nonzero, \(R\se V_\C\), and \(\ep\in\lkv{2}{R^*}\).
It is known that pure spinors, up to multiplication by a constant, are in one to one
correspondence with maximal isotropics \cite{CHEV}.
The maximal isotropics \(L_\mu\) and \(\cjg{L_\mu}=L_{\cjg{\mu}}\) will have trivial
intersection if, and only if, \cite{CHEV}
\[ \muk{\mu}{\cjg{\mu}}\neq0. \]
Here \(\muk{.}{.}\) is the Mukai pairing.  See \cite{CHEV} for details.

For each \(m\in M\) the previous constructions
can be applied to \((\gtsp{m}{M})^\C\).
The bundle formed by these Clifford algebras is called the {\it Clifford bundle}.  In
this context spinors are members of \(\cdf{M}\), and their null spaces are maximal
isotropic subbundles of \(\gtbc{M}\).  The following proposition follows from Theorem 6.4
of \cite{ALEK-XU}.  It is also proven in \cite{GUALT_thesis} for the untwisted case.

\begin{prop}
Generalized complex structures are in one to one correspondence with pure spinor line
bundles \(\Lie\se\lv{\ctbc{M}}\) such that the following hold.
\begin{enumerate*}
  \item  If \(\mu\in\Lie\) then \(\muk{\mu}{\cjg{\mu}}\neq 0\).
  \item  For any local nonzero section \(\mu\in\osn{\Lie\res{U}}\) there exists a local
    section \(X+\xi\in\gvfc{U}\) such that \(d_\Omega(\mu)=(X+\xi)\cdot\mu\).  Here
    \(d_\Omega=d+\Omega\).
\end{enumerate*}
For each point \(m\in M\) the spinor line \(\Lie\res{m}\) is of the form
\eqref{spinorform}
\end{prop}

Let \fct{h}{N}{M} be a twisted
submanifold with a generalized complex structure defined by a spinor line bundle
\(\Lie\se\lv{\ctbc{M}}\).  This spinor line bundle naturally induces a line bundle
in \(\cdf{N}\) given by \(h^*\Lie\).  This induced line bundle could potentially give a
generalized complex structure on \(N\).  We will show that the maximal isotropic defined
by \(h^*\Lie\) is in fact \(\B_h(L_\Lie)\).

\begin{prop}
Let \((M,J,\Omega)\) be a twisted generalized complex manifold and let \(\theta\) be the
pure spinor line that also gives \(J\).  If \fct{h}{N}{M} is a twisted generalized
complex submanifold of \(M\), with generalized complex structure \(J'\), then the spinor
associated to \(J'\) is \(h^*\theta\).
\end{prop}
\begin{proof}
Let \(L\) denote the Dirac structure associated to \(J\).  The spinor line bundle
associated to \(L\) is given by 
\[ \Lie=\set{c(\det(R^o))\wedge\exp(\ep):c\in\C}. \]
Now, 
\[ h^*\big(c\det(R^o)\wedge\exp(\ep)\big)
   =c\det((h^{-1}R)^o)\wedge\exp(h^*\ep). \]
This line bundle is the same as the line bundle associated to \(\B_h(L)\).
\end{proof}

With this proposition it is now a simple matter to give the conditions for a twisted
generalized complex submanifold in terms of spinors.  The involutivity is guaranteed
by \rlem{15}.

\begin{cor}
Let \(M\) be a twisted generalized complex submanifold, with associated spinor line
bundle \(\Lie\).  A twisted submanifold \fct{h}{N}{M} is a twisted generalized complex
submanifold if, and only if, \(h^*\Lie\) is a pure spinor line bundle and
\(\muk{h^*\mu}{h^*\cjg{\mu}}\neq0\) for all \(\mu\in\Lie\).
\end{cor}
%=========================================================================================
\section{Generalized K\"ahler submanifolds}
%=========================================================================================

Finally we will consider submanifolds of generalized K\"ahler structures.  A {\it twisted
generalized K\"ahler structure} on \(M\) is a pair of twisted generalized complex
structures \fct{J_1,J_2}{\gtb{M}}{\gtb{M}} such that:
\begin{enumerate*}
  \item  \(J_1\) and \(J_2\) commute,
  \item  \(\ip{X+\xi}{J_1J_2(Y+\eta)}\), is a positive definite metric.
\end{enumerate*}

The first proposition of this section gives a condition, in terms of the eigenbundles,
for when two complex maps will commute.

\begin{prop}\label{trivial}
Let \(W\) be a real vector space with two maps \fct{\psi_1,\psi_2}{W}{W} such that
\(\psi_1^2=\psi_2^2=-\id\).  Also, let \(L^k_+\) denote the {\pebs} of these maps, and
\(L^k_-\) denote the {\nebs}.  Using this notation, \(\psi_1\) and \(\psi_2\) commute
if, and only if,
\[W_\C=(L^1_+\cap L^2_+)\ds(L^1_+\cap L^2_-)\ds(L^1_-\cap L^2_+) \ds(L^1_-\cap L^2_-).\]
\end{prop}
\begin{proof}
First assume the two maps commute.  Because of this fact, every \(w\in W_\C\) can
be written as
\begin{align*} w
  =& \tfrac{1}{4}(\id-iJ_1)(\id-iJ_2)(w)+\tfrac{1}{4}(\id-iJ_1)(\id+iJ_2)(w) \\
   &\quad+\tfrac{1}{4}(\id+iJ_1)(\id-iJ_2)(w)+\tfrac{1}{4}(\id+iJ_1)(\id+iJ_2)(w) \\
  :=& w^+_++w^+_-+w^-_++w^-_-. 
\end{align*}
It is clear that \(w^\pm_\bullet\in L^1_\pm\), and \(w^\bullet_\pm\in L^2_\pm\).  Now
assume every \(w\in W_\C\) can be written as \(w=w^+_++w^+_-+w^-_++w^-_-\),
where \(w^\pm_\bullet\in L^1_\pm\) and  \(w^\bullet_\pm\in L^2_\pm\).  Now
\((\psi_2\of\psi_1)(w)=\psi_2(iw^+_++iw^+_--iw^-_+-iw^-_+)=-w^+_-+w^+_-+w^-_+-w^-_-\),
and \((\psi_1\of\psi_2)(w)=-w^+_++w^+_-+w^-_+-w^-_-\).
\end{proof}

Let \fct{J_1,J_2}{\gtb{M}}{\gtb{M}} be two commuting bundle maps such that
\(J_1^2=J_2^2=-\id\).  Also, let \(L^k_\pm\) denote the {\pebs} and {\nebs} of \(J_k\).
We also use the notation of \rsect{ind}.  In that \(B=\tb{N}\ds\ctb{M}\res{N}\), and
\fct{s}{B}{\gtb{N}}.
The next two lemmas relate the condition above to our conditions.

\begin{lem}\label{lalg}
The following are equivalent
\begin{enumerate*}
  \item  \(\gtbc{N}=s\big((L^1_+\cap L^2_+)\cap B_\C\big)+
    s\big((L^1_+\cap L^2_-)\cap B_\C\big)+s\big((L^1_-\cap L^2_+)\cap B_\C\big)
    +s\big((L^1_-\cap L^2_-)\cap B_\C\big)\).
  \item  \(\gtb{N}=s(B\cap J_1B\cap J_2B\cap J_1J_2B)\).
  \item  \(B=B\cap J_1B\cap J_2B\cap J_1J_2B+B^\perp\).
  \item  \(J_1B^\perp\cap B\se B^\perp\), \(J_2B^\perp\cap B\se B^\perp\), and
    \(B\cap J_1J_2B^\perp\se B^\perp\).
\end{enumerate*}
\end{lem}
\begin{proof}~
\begin{itemize*}
\item []{\it (1) \(\then\) (2):}
  Every \(v\in\gtbc{N}\) can be written as \(v=s(\vt^+_+)+s(\vt^+_-)+s(\vt^-_+)
  +s(\vt^-_-)\), where \(\vt^+_\pm\in L^1_+\cap L^2_\pm\cap B_\C\) and
  \(\vt^-_\pm\in L^1_-\cap L^2_\pm\cap B_\C\).  Now let
  \(\vt=\vt^+_++\vt^+_-+\vt^-_++\vt^-_-\), and so \(v=s(\vt)\) and \(\vt\in B_\C\).  Now
  \(J_1(\vt)\in B_\C\), \(J_2(\vt)\in B_\C\), and \(J_1J_2(\vt)\in B_\C\).  Thus
  \(\vt\in J_1(B_\C)\), \(\vt\in J_2(B_\C)\), and  \(\vt\in J_1J_2(B_\C)\).  Finally,
  taking the real parts of each of these gives (2).
\item []{\it (2) \(\then\) (1):}  Every \(v\in\gtb{N}\) can be written as
  \(v=s(\vt)\) for some \(\vt\in B\cap J_1B\cap J_2 B\cap J_1J_2B\).  Alternately
  \(v=s(\vt)\) for some \(\vt\in B\) such that \(J_1(\vt)\in B\), \(J_2(\vt)\in B\),
  and \(J_1J_2(\vt)\in B\).  Now we can write
  \begin{multline*}
    \vt=\tfrac{1}{4}\big((\id-iJ_1)\of(\id-iJ_2)(\vt)+(\id-iJ_1)\of(\id+iJ_2)(\vt) \\
       +(\id+iJ_1)\of(\id-J_2)(\vt)+(\id+iJ_1)\of(\id+iJ_2)(\vt)\big).
  \end{multline*}
  By definition each of these components is in the intersection of the eigenbundles, and
  the previous discussion shows that each of these terms is also in \(B_\C\).
\item []{\it (2) \(\iff\) (3):}  We know \(s(B)=\gtbc{N}\), and \(\ker(s)=B^\perp\).
  Thus these two conditions are equivalent.
\item []{\it (3) \(\iff\) (4):}  This equivalence is a fairly straightforward calculation
    \[ B=B^\perp+B\cap J_1B\cap J_2B\cap J_1J_2B \]
  if, and only if,
    \[ B\se B^\perp+J_1B\cap J_2B\cap J_1J_2B \]
  if, and only if,
    \[ B\cap(J_1B^\perp+J_2B^\perp+J_1J_2B^\perp)\se B^\perp \]
  if, and only if,
    \[ B\cap J_1B^\perp+B\cap J_2B^\perp+B\cap J_1J_2B^\perp\se B^\perp \]
  if, and only if,
    \[ B\cap J_1B^\perp\se B^\perp,\, B\cap J_2B^\perp\se B^\perp, \text{ and }
      B\cap J_1J_2B^\perp\se B^\perp. \]
\end{itemize*}
\end{proof}

This last lemma strengthens the conclusions of the first statement in the previous
lemma.

\begin{lem}\label{last}
If \(N\) is a twisted generalized complex submanifold of \((M,J_1)\) and \((M,J_2)\) then
the sums in expression (1), of the previous proposition, are direct.  Also, each of the
components in this expression can be rewritten as
  \[ s\big((L^1_\pm\cap L^2_\pm)\cap B_\C\big)=F^1_\pm\cap F^2_\pm, \]
where \(F^k_\pm=\B_i(L^k_\pm)\).
\end{lem}
\begin{proof}
The fact that \(J_1\) and \(J_2\) descend to generalized complex structures on \(N\)
implies that \(F^k_+\cap F^k_-=\z\), and the sums must be direct.  Now, by definition
\(s(L^k_\pm\cap B_\C)=F^k_\pm\) and it is obvious that
\(s\big((L^1_\pm\cap L^2_\pm)\cap B_\C\big)\se F^1_\pm\cap F^2_\pm\).  To see the other
inclusion, consider \(F^1_+\cap F^2_+\).  This subset will have zero intersection with
\(F^1_-\) and \(F^2_-\), and so it will not intersect with any of the other components.
However, \(F^1_+\cap F^2_+\se\gtbc{N}\), and the fact that \(\gtbc{N}\) is made up of
these four components implies that \(F^1_+\cap F^2_+\se s\big((L^1_\pm\cap L^2_\pm)
\cap B_\C\big)\).
\end{proof}

We are now ready to prove our last theorem, namely that these conditions are guaranteed
to be satisfied by a generalized K\"ahler structure and so our notion of generalized
complex submanifold preserves generalized K\"ahler structures.

\begin{thm}
Let \(N\) be a twisted submanifold of a generalized K\'ahler submanifold \((M,J_1,J_2)\).
If \(N\) is a twisted generalized complex submanifold of \((M,J_1)\) and \((M,J_2)\), then
\((N,J_1',J_2')\) is automatically a twisted generalized K\"ahler manifold.
\end{thm}
\begin{proof}
All that we need to show is \(J_1'J_2'=J_2'J_1'\), and the metric induced by
\(G'=J_1'J_2'\) is positive definite.  We start with the commutativity.  By \rlem{last}
and \rprop{trivial}, if one of the equivalent conditions in \rlem{lalg} is true then
\(J_1\) and \(J_2\) will commute.  Consider condition (4) of this lemma.  By assumption
\(J_1B^\perp\cap B\se B^\perp\) and \(J_2B^\perp\cap B\se B^\perp\).  All that remains
is to show
\(B\cap J_1J_2B^\perp\se B^\perp\).  Take \(v\in B\cap J_1J_2 B^\perp\), so
\(v\in B^\perp\) and \(J_1J_2 v\in B\).  Thus \(\ip{v}{J_1J_2v}=0\).  However, by
assumption this metric is positive definite and so \(v=0\).  Thus
\(B\cap J_1J_2B^\perp\se\z\), and \(B\cap J_1J_2 B^\perp\se B^\perp\) is always true.
It remains to show that \(J_1'J_2'\) defines a positive definite metric.
Take \(v\in\gtb{M}\) and \(\vt\in B\cap J_1B\cap J_2B\cap J_1J_2 B\) such that
\(s(\vt)=v\).  Because \(s\) does not change the inner product
\[ \ip{v}{J'_1J'_2(v)}=\ip{s(\vt)}{J'_1J'_2s(\vt)}=\ip{s(\vt)}{sJ_1J_2(\vt)}
   =\ip{\vt}{J_1J_2(\vt)}, \]
and the positive definiteness of \(J_1J_2\) implies the positive definiteness of
\(J_1'J_2'\).
\end{proof}

%=========================================================================================
%               DOCUMENT END MATTER.
%=========================================================================================
%\bibliographystyle{plainurl}
%\bibliography{/home/jman/latex/maths}

%=========================================================================================
%               DOCUMENT POSTAMBLE.
%=========================================================================================
\end{document}